\magnification=\magstep1
\def\to{\ \longrightarrow\ }

\def\nl{\hfill\break}

\def\hexnumber#1{\ifcase#1 0\or 1\or 2\or 3\or 4\or 5\or 6\or 7\or 8\or
 9\or A\or B\or C\or D\or E\or F\fi}
%
%
\font\twelvemsa=msam10 scaled 1200   
\font\tenmsa=msam10                  
\font\ninemsa=msam9            \font\sevenmsa=msam7
\font\sixmsa=msam6             \font\fivemsa=msam5
%
%
\newfam\msafam                 \textfont\msafam=\tenmsa
\scriptfont\msafam=\sevenmsa   \scriptscriptfont\msafam=\fivemsa
\edef\hexa{\hexnumber\msafam}        
\def\msa{\fam\msafam\tenmsa}         
%
%
\font\twelvemsb=msbm10 scaled 1200   
\font\tenmsb=msbm10                  
\font\ninemsb=msbm9            \font\sevenmsb=msbm7
\font\sixmsb=msbm6             \font\fivemsb=msbm5
%
\newfam\msbfam                 \textfont\msbfam=\tenmsb       
\scriptfont\msbfam=\sevenmsb   \scriptscriptfont\msbfam=\fivemsb
\edef\hexb{\hexnumber\msbfam}        
\def\msb{\fam\msbfam\tenmsb}         
%
%
\font\twelveeufm=eufm10 scaled 1200  
\font\teneufm=eufm10                 
\font\nineeufm=eufm9           \font\seveneufm=eufm7
\font\sixeufm=eufm6            \font\fiveeufm=eufm5
%
\newfam\eufmfam                \textfont\eufmfam=\teneufm
\scriptfont\eufmfam=\seveneufm \scriptscriptfont\eufmfam=\fiveeufm
\edef\hexf{\hexnumber\eufmfam}      
\def\frak{\fam\eufmfam\teneufm}     
%
%
%
\font\twelverm=cmr10 scaled 1200    
\font\ninerm=cmr9                   
\font\sixrm=cmr6   
%
\font\twelvei=cmmi10 scaled 1200    
\font\ninei=cmmi9                   
\font\sixi=cmmi6  
%
\font\twelvesy=cmsy10 scaled 1200   
\font\ninesy=cmsy9                  
\font\sixsy=cmsy6  
%
\font\twelvebf=cmbx10 scaled 1200   
\font\ninebf=cmbx9                  
\font\sixbf=cmbx6  
%
%
\font\twelveit=cmti10 scaled 1200   
\font\nineit=cmti9                  
%
\font\twelvesl=cmsl10 scaled 1200   
\font\ninesl=cmsl9                  
%
\font\twelvett=cmtt10 scaled 1200   
\font\ninett=cmtt9                  
%
%
%
%
\def\small{%
%
%
\textfont0=\ninerm \scriptfont0=\sixrm \scriptscriptfont0=\fiverm
\def\rm{\fam0\ninerm}        
%
%
\textfont1=\ninei \scriptfont1=\sixi \scriptscriptfont1=\fivei
%
%
\textfont2=\ninesy \scriptfont2=\sixsy \scriptscriptfont2=\fivesy
%
%
\textfont3=\tenex \scriptfont3=\tenex \scriptscriptfont3=\tenex
%
%
\textfont\bffam=\ninebf \scriptfont\bffam=\sixbf
\scriptscriptfont\bffam=\fivebf \def\bf{\fam\bffam\ninebf}%
%
%
\textfont\itfam=\nineit \def\it{\fam\itfam\nineit}%
\textfont\slfam=\ninesl \def\sl{\fam\slfam\ninesl}%
\textfont\ttfam=\ninett \def\tt{\fam\ttfam\ninett}%
%
%
%
\textfont\msafam=\ninemsa \scriptfont\msafam=\sixmsa
\scriptscriptfont\msafam=\fivemsa \def\msa{\fam\msafam\ninemsa}%
%
%
\textfont\msbfam=\ninemsb \scriptfont\msbfam=\sixmsb
\scriptscriptfont\msbfam=\fivemsb \def\msb{\fam\msbfam\ninemsb}%
%
%
\textfont\eufmfam=\nineeufm  \scriptfont\eufmfam=\sixeufm
\scriptscriptfont\eufmfam=\fiveeufm \def\frak{\fam\eufmfam\nineeufm}%
%
%
%
\normalbaselineskip=11pt
\setbox\strutbox=\hbox{\vrule height8pt depth3pt width0pt}%
%
%
\normalbaselines\rm}    
%
%
%
%
\def\large{%
\textfont0=\twelverm \scriptfont0=\ninerm \scriptscriptfont0=\sevenrm
\def\rm{\fam0\twelverm}%
\textfont1=\twelvei \scriptfont1=\ninei \scriptscriptfont1=\seveni
\textfont2=\twelvesy \scriptfont2=\ninesy \scriptscriptfont2=\sevensy
\textfont3=\tenex \scriptfont3=\tenex \scriptscriptfont3=\tenex
\textfont\bffam=\twelvebf \scriptfont\bffam=\ninebf
\scriptscriptfont\bffam=\sevenbf \def\bf{\fam\bffam\twelvebf}%
\textfont\itfam=\twelveit \def\it{\fam\itfam\twelveit}%
\textfont\slfam=\twelvesl \def\sl{\fam\slfam\twelvesl}%
\textfont\ttfam=\twelvett \def\tt{\fam\ttfam\twelvett}%
\textfont\msafam=\twelvemsa \scriptfont\msafam=\ninemsa
\scriptscriptfont\msafam=\sevenmsa \def\msa{\fam\msafam\twelvemsa}         
\textfont\msbfam=\twelvemsb \scriptfont\msbfam=\ninemsb
\scriptscriptfont\msbfam=\sevenmsb \def\msb{\fam\msbfam\twelvemsb}         
\textfont\eufmfam=\twelveeufm  \scriptfont\eufmfam=\nineeufm
\scriptscriptfont\eufmfam=\seveneufm \def\frak{\fam\eufmfam\teneufm}
\normalbaselineskip=15pt
\setbox\strutbox=\hbox{\vrule height11pt depth4pt width0pt}%
\normalbaselines\rm}%
%
\def\Bbb{\msb}

%

%
\mathchardef\plussquare="0\hexa01
\mathchardef\nge="3\hexb0B
\mathchardef\maltesecross="0\hexa7A
\mathchardef\del="0\hexf01
%

%

 \input epsf
\overfullrule=0pt
\mathsurround=2pt

\font\npt=cmr9
\font\Bbb=msbm10

\font\secfont=cmbx10

\font\nam=cmr8
\font\aff=cmti8
\font\refe=cmr9

\font\teneufm=eufm10

\mathchardef\square="0\hexa03
\def\qed{\hfill$\square$\par\rm}

\def\boxing#1{\ \lower 3.5pt\vbox{\vskip 3.5pt\hrule \hbox{\strut\vrule
\ #1 \vrule} \hrule} }

\def\down#1{\ \lower 3.5pt\vbox{\vskip 3.5pt \hbox{\strut \ #1 \vrule} \hrule}
}
\def\negdown#1{\ \lower 3.5pt\vbox{\vskip 3.5pt \hbox{\strut  \vrule \ #1
}\hrule} }

\def\Hom{\mathop{\rm Hom}}

\hsize=6.3 truein
\vsize=9 truein

\baselineskip=13 pt
\parskip=\baselineskip
 1

\parindent=0pt

\def\Z{\hbox{\Bbb Z}}

\def\x{\hbox{\bf x}}
\def\y{\hbox{\bf y}}

\def\M{{\cal M}}







\newif \iftitlepage \titlepagetrue

\def\diagram{\global\advance\diagramnumber by 1
$$\epsfbox{turfig.\number\diagramnumber}$$}
\def\ddiagram{\global\advance\diagramnumber by 1
\epsfbox{turfig.\number\diagramnumber}}

\newcount\diagramnumber
\diagramnumber=0

\newcount\secnum \secnum=0
\newcount\subsecnum
\newcount\defnum
\def\section#1{
               \vskip 10 pt
               \advance\secnum by 1 \subsecnum=0
               \leftline{\secfont \the\secnum \quad#1}
               }

\def\subsection#1{
               \vskip 10 pt
               \advance\subsecnum by 1
               \defnum=1
               \leftline{\secfont \the\secnum.\the\subsecnum\ \quad #1}
               }

\def\definition{
               \advance\defnum by 1
               \bf Definition
\the\secnum .\the\defnum \rm \
               }

\def\lemma{
               \advance\defnum by 1
               \par\bf Lemma  \the\secnum
.\the\defnum \sl \ \par \rm
               }

\def\theorem{
               \advance\defnum by 1
               \par\bf Theorem  \the\secnum
.\the\defnum \rm \ \par
              }

\def\cite#1{\secfont [#1]\rm}

\vglue 20 pt

\centerline{\secfont WEYL ALGEBRAS and KNOTS}

\medskip

\centerline{\nam ROGER FENN${}^1$, VLADIMIR TURAEV${}^2$}
\centerline{\aff ${}^1$School of Mathematical Sciences, University of Sussex}
\centerline{\aff Falmer, Brighton, BN1 9RH, England}
\centerline{\aff e-mail address: rogerf@sussex.ac.uk}
\centerline{\aff ${}^2$Institut de Recherche Math\'ematique Avanc\'ee, CNRS --}
\centerline{\aff Universit\'e Louis Pasteur, 7 rue Ren\'e Descartes, }
\centerline{\aff 67084 Strasbourg Cedex France}
\centerline{\aff e-mail address: turaev@math.u-strasbg.fr }

\baselineskip=10 pt
\parskip=0 pt
\bigskip
\centerline{\nam ABSTRACT}
\leftskip=0.25 in
\rightskip=0.25in

{\nam In this paper we push forward results on the invariant
${\cal F}$-module of a virtual knot investigated by the first
named author where ${\cal F}$ is the algebra with two invertible
generators $A,B$ and one relation $A^{-1}B^{-1}AB-B^{-1}AB=
BA^{-1}B^{-1}A-A$. For flat knots and links the two sides of the 
relation equation are put equal to unity and the algebra becomes
the Weyl algebra. If this is perturbed and the two sides of the 
relation equation are put equal to a general element, $q$, of the ground
ring, then the resulting module lays claim to be the correct
generalization of the Alexander module. Many finite dimensional
representations are given together with calculations.
}

\leftskip=0 in
\rightskip=0 in
\baselineskip=13 pt
\parskip=\baselineskip

\parskip=\baselineskip
\section{Introduction}
The first named author \cite{F} and others recently introduced a
general method associating a module to an arbitrary oriented link
diagram on an oriented surface. The module is generated by the arcs
obtained by splitting the diagram at its crossings. The relations are
associated with the crossings and depend on a choice of two invertible
elements $A,B$ of an associative algebra such that
$$A^{-1}B^{-1}AB-B^{-1}AB=BA^{-1}B^{-1}A-A$$
This, {\it fundamental relation}, ensures that the module is invariant under
the
Reidemester moves on the diagram  and provides thus an invariant of oriented
links in surfaces crossed with an interval.

This construction generalizes the classical Alexander module of a link
in $S^3$ as well as its natural extensions to surface links. Note
however that the classical construction uses commuting $A,B$ which
conceals non-commutative ramifications.

In the present paper we show that the fundamental equation has a
natural family of solutions arising from the so-called quantum Weyl
algebras. For physicists this is a quantum version of the Heisenberg
algebra of a harmonic oscillator. Given an element $q$ of a
(commutative) ground ring $K$, we define the {\it quantum Weyl
algebra} $W_q$ to be the $K$-algebra generated by four elements $u,
u^{-1}, v, v^{-1}$ subject to the relations
$uu^{-1}=u^{-1}u=vv^{-1}=v^{-1}v=1$ and $uv-qvu =1$.

We show that $A=v^{-1}u^{-1}\in W_q$ and $ B=u\in W_q$ satisfy the
fundamental relation.  Applying Fenn's method to these $A,B$, we
associate with every link diagram $D$ a $W_q$-module $\M(D)$. It has a
square presentation matrix so that one is tempted to take the
determinant (or subdeterminants) as in the Alexander theory. However,
the algebra $W_q$ is non-commutative and the determinants do not make
sense. One solution is to plug in matrix representations of $W_q$ over
$K$ and take the determinants of the resulting matrices, see for
example, \cite{As}. This gives interesting link invariants satisfying
the same skein relation as the Alexander-Conway polynomial of links in
$S^3$ and dependant on the choice of a matrix representation of
$W_q$. As expected, these invariants are extensions of the known
Alexander invariants of links in $S^3$.

Parallel constructions work for closed curves on oriented surfaces and
produce homotopy invariants of such curves. Here we involve an
extension of the classical Weyl algebra, obtained by putting $q=1$ in
$W_q$.

The plan of the paper is as follows. In Sections 2 and 3 we discuss
the extension of classical Weyl algebra and use it to produce
invariants of closed curves on surfaces. In Sections 4 and 5 we
discuss the quantum Weyl algebra and use it to produce link
invariants. In Section 6 we give proofs of several claims made in
Sections 3 and 5. In Section 7 we consider representations of the
virtual braid groups arising from quantum Weyl algebras. In Section 8
some calculations are given. One of the examples shows that the shadow
of the Kishino knot is non-trivial verifying a result in
\cite{Kad}.

The authors are grateful to Vladimir Bavula for a useful discussion of
Weyl algebras and their quantum counterpoints used in section 4.

\section{The Weyl Algebra}
In this section we consider the Weyl algebra and variants. A good
reference for the details of this section is Cohn's book \cite{Co}.

Let $K$ be a commutative ring and $W^0$ be the $K$-algebra generated
by $u, v$ and with relation
$$uv-vu=1.$$ Then $W^0$ is called the {\it Weyl algebra} on $u, v$
over $K$.  If $K$ is a field of characteristic 0 then $W^0$ is a
simple ring, that is all two sided ideals are trivial, see
\cite{Co}, pp 362-363.

We now extend $W^0$ so that $u, v$ are invertible. Let $W$ be the {\it
extended Weyl algebra} defined as the quotient of the  $K$-algebra  
generated by $u^{\pm 1}, v^{\pm 1}$ by the ideal generated by 
$uv-vu-1$. Although we shall not formally need it, note
that the natural algebra homomorphism $W^0\to W$, sending $u$ to $u$
and $v$ to $v$, is injective. In section 4 we will give a proof of a more
general statement.

Let $M_n(K)$ denote the algebra of $n\times n$ matrices with entries
in $K$.  We would like to represent $W$ by matrices in $M_n(K)$. Note
that, given $v$, the relation $uv-vu=1$ is affine in $u$ and so any
solution can be written $u=u_P+u_H$ where $u_P$ is any particular
solution and $u_H$ commutes with $v$, say a polynomial in $v$.

\theorem{\sl The extended Weyl algebra $W$ over a field of characteristic
zero $K$  has no non-trivial representations in $M_n(K)$.}

{\bf Proof} Since $W^0$ is simple, i.e., all two sided ideals are
trivial, any representation of $W^0$ is either trivial or
faithful. But $W^0$ has infinite dimension over $K$ and so any finite
dimensional representation collapses on $W^0$ and hence on $W$.\qed

An obvious question is whether $W$ is simple over a field of
characteristc zero.  Note that $W$ is not simple if the field has
characteristic a prime $p$ since $K[u^p,v^p]$ is central.

Traditionally $W$ acts on the algebra $C^\infty({\bf C}-\{0\})$ by
$$u(f)=f'+f,\quad v(f)=xf$$
where $f\in C^\infty({\bf C}-\{0\})$.

For a finite dimensional representation, consider the truncated polynomial ring
$$R=K[x]/(x^n=0)$$
where $K$ is a field of characteristic dividing $n$.
Let $I=i_0+i_1x+\cdots$ and $J=j_0+j_1x+\cdots$ be elements of $R$. If
$i_0\ne 0$ and $i_1\ne 0$ then $I$ and its derivative $I'$ are units.
In fact $(I')^{-1}=k_0+k_1x+k_2x^2+\cdots$ where
$$k_0=i_1^{-1},\ k_1=-2i_2i_1^{-2},\ k_2=(4i_2^2-3i_1i_3)i_1^{-3},\ldots$$
Further coefficients can be determined from the difference equation
$$i_1k_r+2i_2k_{r-1}+\cdots+(r+1)i_{r+1}k_0=0$$
Define the $K$-linear operators $u,v:R\to R$ by
$$u(f)={f'\over I'}+Jf,\quad v(f)=If.$$
Then it is easily seen that $  v$ is invertible and $uv-vu=1$.
The matrices of $u, v$ with respect to the basis $\{1, x, x^2, \ldots,
x^{n-1}\}$ are
$$u=\pmatrix{
j_0&j_1&j_2&\ldots&j_{n-2}&j_{n-1}\cr
k_0&j_0+k_1&j_1+k_2&\ldots&j_{n-3}+k_{n-2}&j_{n-2}+k_{n-1}\cr
0&2k_0&j_0+2k_1&\ldots&j_{n-4}+2k_{n-3}&j_{n-3}+2k_{n-2}\cr
\vdots&\vdots&\vdots&\ddots&\vdots&\vdots\cr
0&0&0&\ldots&(n-1)k_0&j_0+(n-1)k_1\cr}$$
$$v=\pmatrix{
i_0&i_1&\ldots&i_{n-1}\cr
0&i_0&\ldots&i_{n-2}\cr
0&0&\ldots&i_{n-3}\cr
\vdots&\vdots&\ddots&\vdots\cr
0&0&\ldots&i_0\cr}$$
Then $u$ is invertible provided its determinant is
non-zero. This   imposes a polynomial condition
on the variables, $i_0, i_1,\ldots,j_0, j_1,\ldots$.

If $K$
has characteristic $p$ dividing $n$, then
$$u=\pmatrix{
x&a_1&0&\ldots&0\cr
0&x&a_2&\ldots&0\cr
\vdots&\vdots&\vdots&\ddots&\vdots\cr
0&0&0&\ldots&a_{n-1}\cr
0&0&0&\ldots&x\cr}\quad
v=\pmatrix{
y&0&\ldots&0&0\cr
1/a_1&y&\ldots&0&0\cr
0&2/a_2&\ldots&0&0\cr
\vdots&\vdots&\ddots&\vdots&\vdots\cr
0&0&\ldots&(n-1)/a_{n-1}&y\cr}$$
is a non-trivial representation  of $W$ in $M_n(K)$, with $n+1$
parameters $x,y ,  a_1,..., a_{n-1}\in K$ where $x\neq 0$, $ y \neq 0$.

\section{$W$-Modules and Flat Links}
By a {\it system of loops on a surface}, we mean a mapping from a disjoint
union of a finite number of oriented circles   into a
compact oriented surface   with empty boundary.  Two systems of
loops on surfaces are {\it stably equivalent} if they can be related
by a finite sequence of the following operations: (i) homotopy of
loops on the surface; (ii) composing   with an orientation
preserving homeomorphism of surfaces; (iii) attaching a 1-handle to
the ambient surface away from the loops or removing such a handle.

We   call the stable equivalence classes of systems of loops on surfaces {\it
flat links}. When there is only one loop, we can  speak  of {\it
flat knots}. Elsewhere they are called flat virtuals, \cite{Ka} and
virtual strings, \cite{Tu}.


We now define a $W-$module, for each flat link $L$.  Represent $L$ by a system
of loops on a surface lying   in general position. This system has  a finite
number, $n$, of  intersection   points, where the loops cross (or self-cross)
transversely. Let $m$ be the number of loops in the system having no
self-intersection points and missing the other loops. The {\it arcs} are
the $N=2n+m$ components of our  system of loops with the   intersection points
removed. Note that arcs do not pass through   crossing
points. Label the arcs  $x_1,\ldots, x_{N}$ in an arbitrary way.
These labels will be the generators of the module.

Pick an associative $K$-algebra with unit $R$ and fix    elements $A,B,C,D \in
R$.  Each
self-intersection point contributes 2 linear relations between the symbols
$x_1,\ldots, x_{N}$ as indicated by the
following diagram where the input arcs arriving at the crossing are labelled
$x_i, x_j$, the
output arcs leaving the crossing are labelled $x_k, x_l$ and  we assume that
the surface is oriented
counterclockwise.
\diagram

Quotienting the free $R$-module with generators $x_1,\ldots, x_{N}$
by the $2n$ relations derived from the crossings, we obtain an $R$-module
denoted  $\M_{A, B, C, D}(L)$.

\theorem {\sl Let $R=W$ and   $A=v^{-1}u^{-1}$, $B=u$,
$C=uvu^{-1}v^{-1}u^{-1}vuv^{-1}u^{-1}$, $D=-u^{-1}v^{-1}$ (where
$uv-vu=1$). Then the
$W-$module $\M_{A, B, C, D}(L)$ is a stable equivalence invariant of $L$.}

This theorem will be proven in Section 6.

Let  $\M(L)$ denote the $W-$module defined in Theorem 3.2. Workable invariants
can be derived from
this module as follows.
If $\M'$ is any $W-$module then $\Hom_W(\M(L),\M')$ is a finitely generated
$K$-module. If $K$ is a field, we can set
$$d_{\M'}(L)=\dim_K{\Hom}_W(\M(L),\M').$$
Then $d_{\M'}(L)$ is an integer invariant of $L$. For the trivial $m$-component
flat link, $\M(L)=W^m$ and  $d_{\M'}=(\dim_K(\M'))^m$.

A presentation of the module $\M(L)$ as above is determined by an
$N\times  N$   matrix $M$ with entries in $W$  (actually the entries are $A, B,
C,
D, 0,-1$).   Any  matrix representation  $\rho:W\to M_k(K)$ transforms   $M$
into a $2Nk\times 2Nk$ square matrix. Its
determinant $\Delta_0$ is an invariant of $L$ up to multiplication by
powers of the determinant of $\rho(B)$, see \cite{F}. We can also
consider the ideal $I_r$ in $K$ generated by the codimension $r$
subdeterminants for integer $r>0$.  For suitable $K$, this has a
greatest common divisor $\Delta_r(L)$ which is   an invariant of $L$
up to multiplication by units.

\section{The Quantum Weyl Algebra}
Let $W_q^0$ be the  algebra over the commutative ring $K$ generated by
$u, v$ and with relation $uv-qvu=1$ where $q$ is an invertible element of $K$.
This algebra is called the
{\it $q$-oscillator algebra} on $u, v$  in \cite{DP}. The
variable $h$ used in \cite{DP} is equal to $1$ here.

We define the {\it   quantum 
Weyl algebra} $W_q$ to be the quotient of the $K$-algebra in $u^{\pm
1}, v^{\pm 1}$, by the ideal generated by $uv-qvu-1$.

\lemma{\sl  The algebra homomorphism $\alpha:W_q^0\to W_q$, sending
$u$ to $u$ and $v$ to $v$, is injective.}

{\bf Proof} Recall first the definition of the ring of non-commuting
polynomials, $R[x; \sigma]$, where $R$ is a ring and $\sigma$ is a ring
automorphism of $R$. The ring $R[x;
\sigma]$ is obtained from the free algebra on $R$ with an added
generator $x$ by imposing the condition
 $xr=\sigma(r) x$ for all $ r\in R$.
Any element of  $R[x; \sigma]$ can be written uniquely as a
``finite polynomial'' in $x$, namely $r_0+r_1x+r_2x^2+\cdots+r_nx^n$ where
all $r_i$ lie in $R$. Similarly, we define the ring of non-commuting 
Laurent polynomials, $R[ x^{\pm 1}; \sigma]$, where any element is a Laurent
polynomial $r_{-m}x^{-m}+\cdots+r_nx^n$. Of course if $\sigma$ is the
identity, then these definitions give the usual ring of polynomials, $R[x]$,
and the ring of Laurent polynomials, $R[  x^{\pm 1}]$.

Let $R=K(h)$ be the field of rational functions on one variable $h$
 with coefficients in $K$.  In other words, $R$ is the field of
 fractions of the commutative ring of polynomials $K[h]$.  Let
 $\sigma$ be the ring automorphism of $R$ sending $h$ to
 $q^{-1}(h-1)$. Thus, $\sigma$ sends an arbitrary rational function
 $f(h)\in R$ to $f(q^{-1}(h-1))\in R$. Consider the ring of
 non-commuting Laurent polynomials, $U=R[x^{\pm 1};
\sigma]$. In $U$ we have the equalities 
$$(hx^{-1}) x - q x (hx^{-1})=h- q \sigma(h) x x^{-1}= h-(h-1)=1.$$
This implies the existence of a homomorphism of $K$-algebras $\beta:
W_q^0 \to U$ such that $\beta (u)=hx^{-1} $ and $\beta (v)=x$.  The
same equalities and the fact that $hx^{-1} $ and $x$ are invertible in
$U$ imply that there is a homomorphism of $K$-algebras $\gamma: W_q
\to U$ sending $u$ to $hx^{-1} $ and $v$ to $x$. It is clear that
$\beta =\gamma \alpha$. It is easy to see that $\beta$ is injective.
Indeed, if a polynomial $\sum_{m,n\geq 0} k_{m,n} u^m v^n$ with
$k_{m,n}\in K$ lies in the kernel of $\beta$, then $\sum_{m,n\geq 0}
k_{m,n} h^m x^{n-m}=0$ in $U$. For any integer $s$, the monomial $x^s$
appears here with coefficient $\sum_{m \geq 0} k_{m,m+s} h^m $ which
therefore must be $0$ in $R=K(h)$.  Hence $k_{m,n}$ are all equal to
$0$.  The injectivity of $\beta$ and the equality $\beta =\gamma
\alpha$ imply that $\alpha$ is injective.
\qed

Of course if $q=1$ then $W^0_q$ and $W_q$ reduce to $W^0$ and $W$
considered earlier and the above lemma encompasses them also. From now
on we shall assume that $1-q$ is invertible.

We now study finite dimensional representations of $W_q$.

\lemma{\sl Let $u,v$ be invertible $n\times n$ matrices  over $K$ satisfying
$uv-qvu=1$. If the linear map  $K^n\to K^n$ defined
by $x\mapsto   vx$ has no invariant subspaces other than $0$ and $K^n$,
then either   $u=(1-q)^{-1}v^{-1}$  or $q$ is an $n$-th
root of unity.}

{\bf Proof}
Set
$$u=u_H+{1\over 1-q}v^{-1}$$
Then $u$ is defined by $u_H$ and conversely. Moreover $u_H$
$q$-commutes with $v$,
$$u_Hv=qvu_H.$$
Let $X$ be the kernel of the linear map  $K^n\to K^n$ defined
by   $x\mapsto u_H x$. Then $X$ is an
invariant subspace of $v$ and so is either $K^n$ or $\{0\}$. In the
first case $u=(1-q)^{-1}v^{-1}$. In the second case take determinants. This
shows that $q^n=1$. \qed

In view of this lemma, we will restrict our attention to triangular
matrices. Suppose that $u, v$ are upper triangular with
diagonal elements $a_1, a_2, \ldots, a_n$ and $b_1, b_2,\ldots, b_n$
repectively. Then it is easy to see  that the $b$'s are related to the
$a$'s by
$$a_ib_i={1\over 1-q}.$$
Moreover if $a_i\ne qa_{i+1}, b_{i+1}\ne qb_i$ for all $i$ then
$u, v$ commute. An example where this doesn't happen is
$$u=\pmatrix{
q^{n-1}a&b^{n-2}d&0&\ldots&0\cr
0&q^{n-2}a&b^{n-3}d&\ldots&0\cr
\vdots&\vdots&\vdots&\ddots&\vdots\cr
0&0&0&\ldots&d\cr
0&0&0&\ldots&a\cr},\quad
v=\pmatrix{c&e&0&\ldots&0\cr
0&qc&qb^{-1}e&\ldots&0\cr
\vdots&\vdots&\vdots&\ddots&\vdots\cr
0&0&0&\ldots&(qb^{-1})^{n-2}e\cr
0&0&0&\ldots&q^{n-1}c\cr}
$$
where $c=1/(aq^{n-1}(1-q))$.

Let us look for representations of the form
$$u=\pmatrix{
a_1&0&0&\ldots&0\cr
b_1&a_2&0&\ldots&0\cr
0&b_2&a_3&\ldots&0\cr
\vdots&\vdots&\vdots&\ddots&\vdots\cr
0&0&0&\ldots &a_n\cr},\quad
v=\pmatrix{
c_1&d_1&0&\ldots&0\cr
0&c_2&d_2&\ldots&0\cr
0&0&c_3&\ldots&0\cr
\vdots&\vdots&\vdots&\ddots&\vdots\cr
0&0&0&\ldots&c_{n}\cr}$$
These will satify $uv-qvu=1$ provided
$a_i=q^{n-i}a$, $c_i=q^{n-i}c$ for some $a, c$ and
$$\eqalign{b_id_i&=q^{n-2i}+q^{n-2i+1}+\cdots+q^{n-i-1}-q^{-i+1}-q^{-i+2}-\cdots-q^{-1},\cr
&\qquad\ i=1,2,\ldots,n-1\cr}$$
giving a representation of $W_q$ with $n+1$ parameters.

\section{Virtual Links}

We will consider links in 3-manifolds obtained by thickening of a compact
oriented surface, $\Sigma$,   with empty boundary. These links are represented
by
smooth embeddings  of a disjoint union of a finite number of oriented circles
into $  \Sigma\times [0,1]$.  Under the projection $\Sigma\times
[0,1]\to\Sigma$ a link
defines a link diagram   on $\Sigma$ which we will assume in
general position with only transverse double points. The over and
under arcs can be distinguished in the usual way. This establishes a bijection
between ambient isotopy classes of links in $\Sigma\times [0,1]$ and
link diagrams on $\Sigma$ modulo the Reidemeister moves.

Two  links are said to be {\it stably equivalent} if they can be
related by a finite sequence of the following operations: (i)
Reidemeister moves on a link diagram on the surface; (ii)
transferring a diagram with an orientation preserving homeomorphism of
surfaces; (iii) attaching a 1-handle to the ambient surface away from
the diagram or removing such a handle.  The stable equivalence class
of a link $L$ is a {\it surface link} or a {\it virtual link}, see
\cite{Ka},
\cite{Kam}. We will shorten this to just link.
 There is a natural map from links to flat links defined by
forgetting the over and under information on a link diagram.

In   generalization of the method used for flat links we define a
$W_q-$module, $\M_q(L)$, for each link $L$. It is called the {\it
quantum Weyl} module of $L$ and is defined by
generators and relations using a   diagram of $L$.  Suppose that the diagram
  has $n$ double points and $m$ simple closed loops  disjoint from  other
loops. The {\it arcs} are
the $N=2n+m$ components of the diagram with the   intersection points
removed.   Let the arcs be labelled $x_1,\ldots, x_{N}$.  These
labels will be the generators of the module.

 Pick an associative $K$-algebra with unit $R$ and fix      elements $A,B,C,D
\in R$.  Each
self-intersection point contributes 2 relations as indicated by the
following diagram.
\diagram
We get an $R$-module as we did for flat links which again we denote by $\M_{A,
B, C, D}(L)$.
\theorem{\sl Let $R=W_q$  and $A=v^{-1}u^{-1}$, $B=u$,
$C=quvu^{-1}v^{-1}u^{-1}v^{-1}u^{-1}vu$, $D=1-q-u^{-1}v^{-1}$ (where
$uv-qvu=1$). Then the
$W_q$-module $\M_{A, B, C, D}(L)$ is a stable equivalence invariant of $L$.}

The proof will be given in Section 6.  As in the case of flat links we
may define rank invariants $d_{\M'}(L)$, ideal invariants $I_r^{(q)}$
and determinental invariants $\Delta_0$ and $\Delta_r$.

\section{Proof of Theorems 3.2 and 5.5}

 We say that elements $A, B  $ of an associative $K$-algebra with unit $R$
satisfy the {\it fundamental} relation if they are invertible and
$$\eqalignno{A^{-1}B^{-1}AB-B^{-1}AB&=BA^{-1}B^{-1}A-A.&(1)}$$
Recall the following theorem of Fenn \cite{F, BF, BuF}.

\theorem{\sl If $A, B\in R$ satisfy  the fundamental relation,
$C=A^{-1}B^{-1}A(1-A) $ is invertible in $R$, and $D=1-A^{-1}B^{-1}AB $, then
for any link $L$, the module $\M_{A,B,C,D}(L)$ is a stable equivalence
invariant of $L$.}

To deduce Theorem 5.5 from Theorem 6.6  we need the following lemma.

\lemma{\sl If $q\in K$ and $A, B$ are invertible elements of $R$   satisfying
the relation
$$\eqalignno{B&=BA^{-1}-qA^{-1}B,&(2)}$$
then $A, B$ satisfy the fundamental equation (1).}

{\bf Proof } We have
$$qB^{-1}=A^{-1}B^{-1}A-B^{-1}A $$ as can be seen by
multiplying (2) on the right by $B^{-1}A $ and then multiplying on
the left by $B^{-1}$.  Then the requirement of the fundamental equation is that
$qB^{-1}$ commutes with $B$ which holds since $q$ is an element of the ground
ring. \qed

We apply this lemma to $R=W_q$ and   $A=v^{-1}u^{-1}, B=u$. It is easy to check
that
$A, B$ satisfy (2) and therefore $A,B$ satisfy the fundamental relation. We
have
 $$C=A^{-1}B^{-1}A(1-A)=uv u^{-1}v^{-1}u^{-1} (1- v^{-1}u^{-1})=
quv u^{-1}v^{-1}u^{-1}v^{-1}u^{-1} vu$$
where we use the equality
$1-v^{-1}u^{-1}= q v^{-1}u^{-1} vu$ obtained from $uv-qvu=1$ via
multiplication  by  $v^{-1}u^{-1}$. Note that $C$ is invertible in $W_q$.
Similarly,
$$D=1-A^{-1}B^{-1}AB=1- uv u^{-1}v^{-1}=1-q- u^{-1}v^{-1}.$$
Therefore    Theorem 5.5 follows from Theorem 6.6.

To prove Theorem  3.2   it is enough to substitute $q=1$ in  Theorem 5.5 and to
observe that for our choice of $A,B, C, D \in W$, the matrix $$S=\pmatrix{A & B
\cr C & D\cr}\in M_2(W) $$
is equal to its inverse: $S=S^{-1}$. The latter is a direct consequence of
Claim (c) of the following lemma (Claims (a) and (b) of this lemma will be used
in the next section).

\lemma{\sl Suppose $A, B\in R$ satisfy (1) and $C, D$ are defined by
$$\matrix{%
C=A^{-1}B^{-1}A(1-A),&
D=1-A^{-1}B^{-1}AB\cr}.$$
(a) Let $$S=\pmatrix{A & B \cr C & D\cr}\in M_2(R).$$
Then $S$ is invertible if and only if $C$ is invertible.

(b) Put
$$S_1=\pmatrix{A & B & 0\cr
C & D & 0\cr 0 & 0 & 1\cr}\in M_3(R),\ S_2=\pmatrix{1 & 0 & 0\cr
0 & A & B \cr 0 & C & D\cr}\in M_3(R).$$
Then $$\eqalignno{S_1S_2S_1&=S_2S_1S_2.&(4)}$$

(c) Suppose that $q=(1 - A)A^{-1}B^{-1}AB$ is an element of the ground ring.
Then $S^2=(1-q)S+q$. }

{\bf Proof }

We can write $S$ as a product of elementary matrices
$$\pmatrix{%
A & B \cr C & D\cr}=
\pmatrix{%
A & 0 \cr 0 & 1}
\pmatrix{%
1 & 0 \cr C & 1\cr}
\pmatrix{%
1 & 0 \cr 0 & 1-A^{-1}\cr}
\pmatrix{%
1 & A^{-1}B \cr 0 & 1}.
$$
If $C$ is invertible then so is $1-A$ and hence each matrix in the product
is invertible. This proves (a).

The proof of (b) follows by basic manipulations.

To prove (c) observe that conjugating the equality $q=(1 - A)A^{-1}B^{-1}AB$ by
$B$
we obtain $q=BA^{-1}B^{-1}A -A$. Multiplying  on the right by $A^{-1}$ we
obtain that
$BA^{-1}B^{-1}=1+qA^{-1}$. The latter formula will be used in the computation
of
$$S^2=\pmatrix{A^2+BC & AB+BD\cr CA+DC& CB+D^2\cr}.$$
Substituting the values of $C, D$ we get
$$\eqalignno{
A^2+BC&=A^2+(BA^{-1}B^{-1})A(1-A)\cr
&=q+(1-q)A,\cr
AB+BD&=AB+B-(BA^{-1}B^{-1})AB\cr
&=(1-q)B,\cr
CA+DC&=A^{-1}B^{-1}A^2(1-A)+(1-A^{-1}B^{-1}AB)A^{-1}B^{-1}A(1-A)\cr
&=A^{-1}B^{-1}A(1-A^2-(BA^{-1}B^{-1})A(1-A))\cr
&=A^{-1}B^{-1}A(1-q)(1-A)\cr
&=(1-q)C,\cr
CB+D^2&=A^{-1}B^{-1}A(B-AB+(BA^{-1}B^{-1})AB)+2D-1\cr
&=A^{-1}B^{-1}A(1+q)B+2D-1\cr
&=(1+q)A^{-1}B^{-1}AB+2D-1\cr
&=q+(1-q)D.\cr}$$
As an aid to the reader we have put brackets where the substitution
$BA^{-1}B^{-1}=1+qA^{-1}$ takes place.
\qed

 Note: the quadratic equation in $S$ given by (c), implies that we
have a representation of the Hecke algebra, $H_n$, for each $n$,
\cite{H}. The equation can also be written
$$q^{-1/2}S-q^{1/2}S^{-1}=q^{-1/2}-q^{1/2}.$$

\section{Representations of the Braid Group}
In this section we look at some representations of the braid group,
$B_n$, an extension of the braid group, $VB_n$, and a quotient of this
extension, $FB_n$. These representations are defined by the work in
the previous section. In the case of the braid group all the
representations are equivalent to the Burau representation although
this is certainly not the case for the two other groups.

Let $n$ be a positive integer. The {\it braid group} $B_n$ has
generators $\sigma_1, \ \sigma_2, \ \ldots, \sigma_{n-1}$ and relations
$$\matrix{\hfill \sigma_i \sigma_j&= \sigma_j \sigma_i\hfill &|i-j|>1 \cr
\hfill\sigma_i \sigma_{i+1} \sigma_i&=
\sigma_{i+1}\sigma_i\sigma_{i+1}\hfill & i=1,\ldots, n-1\cr} $$

The {\it virtual braid group}, $VB_n$, is an extension of $B_n$ with
new generators\break $\tau_1, \ \tau_2, \ \ldots, \ \tau_{n-1}$ and two
sorts of extra relations.\hfill\break
Permutation group relations:
$$\matrix{
\hfill {\tau_i}^2&=1 \hfill& \cr
\hfill \tau_i \tau_j&= \tau_j \tau_i \hfill &|i-j|>1 \cr
\hfill \tau_i \tau_{i+1} \tau_i&=\tau_{i+1} \tau_i \tau_{i+1}\hfill \cr } $$
Mixed relations:
$$\matrix{ \hfill \sigma_i \tau_j&= \tau_j \sigma_i &|i-j|>1 \cr
\hfill\sigma_i \tau_{i+1} \tau_i&=\tau_{i+1} \tau_i \sigma_{i+1}\hfill
& i=1,\ldots, n-1\cr } $$
The {\it flat braid group}, $FB_n$, is the quotient of $VB_n$ by the relations
$\sigma_i^2=1$ for all $i$.

Let $S:X^2\to X^2$ be a permutation of the cartesian square of a set $X$.
In \cite{FJK} such an $S$ is called a {\it switch} if
$$(S\times id)(id\times S)(S\times id)= (id\times S)(S\times id)(id\times S).$$

Examples of switches are the identity and the {\it twist}, $T$, defined by
$T(a,b)=(b,a)$.

A binary operation, $(a,b)\to a^b$, is called {\it invertible on the right}
if there exists another binary operation,  $(a,b)\to a^{b^{-1}}$ such that
$$ a^{bb^{-1}}= a^{b^{-1}b}=a$$
is always true. For example racks or quandles, see \cite{FR}, are
invertible.

A switch, $S$, defines two binary operations by the formula
$$S(a,b)=(b_a,a^b).$$
A switch is called a {\it biquandle} if both operations are invertible
and
$$a^{a^{-1}}=a_{a^{a^{-1}}}$$
for all $a\in X$. Note that the original definition included the
extra condition
$$a_{a^{-1}}=a^{a_{a^{-1}}}$$
but this has been shown to be unnecesary by \cite{S}.

\theorem{\sl A linear switch defined by a $2\times2$ matrix is a biquandle}

{\bf Proof} See \cite{FJK} \qed

If $S:X^2\to X^2$ is a switch and $n$ is a positive integer, define
the permutation, $S_i$, of $X^n$ by
$$S_i=id^{i-1} \times S \times id^{n-i-1}$$ where $id:X\to X$ is the
identity.  Then $S$ defines a homomorphism $r_S:B_n\to P_n(X)$ where
$P_n(X)$ is the group of permutations of $X^n$.

Let $b(t)=\pmatrix{0&1\cr t&1-t\cr}$ denote the Burau matrix with
parameter $t$. Then $b(t)$ is a switch.

\theorem{\sl If $S=\pmatrix{A&B\cr C&D\cr}$ is any linear switch, then
the homomorphism $r_S$ is equivalent to $r_{b(t)}$ where
$t=(1-A)(1-D)$.}

{\bf Proof} See \cite{F}. \qed

We now extend the representation $r_S$ to $VB_n$ by sending the
generator $\tau_i$ to $T_i$, where $T$ is the twist. We will continue to call
the resulting
representation $r_S$.

Any   link  $L$ is the closure of a
virtual braid $\beta\in VB_n$ for some $n$. Let $S=\pmatrix{A&B\cr
C&D\cr}$ be a linear switch where $A, B, C, D$ are elements of an
algebra $R$. Let $\M_S(L)$ denote the left $R$-module with $2n\times 2n$
presentation matrix $r_S(\beta)-1$.

For example, suppose that $S=\pmatrix{1-BC&B\cr C&0\cr}$ then the
module $\M_S(L)$ is the Sawollek module of $L$ with a change of variable, see
\cite{Sa}. This becomes the Alexander module of $L$ if $C=1$.

\theorem{\sl For any linear switch $S$, the module $\M_S(L)$ is a stable
equivalence invariant of $L$. If the link $L$ is
classical then $\M_S(L)$ is equivalent to the Alexander module of $L$.  So in
particular $\Delta_0=0$ and $\Delta_1$ is the Alexander polynomial
with variable $t=(1-A)(1-D)$.}

{\bf Proof} For the invariance of the module see \cite{F} or
\cite{FJK}. For classical links the braid $\beta$ can be chosen in
$B_n$ and then the representation is equivalent to the Burau
representation and this defines the Alexander module. \qed

More generally we have the following,
\theorem{\sl Let $\M_q(L)$ be the quantum Weyl module of a link $L$. Then
if the generators $u, v$ of the Weyl algebra commute, the module
becomes the Sawollek module $\M_S(L)$ for $S=\pmatrix{1-q&u\cr q/u&0\cr}$}

Suppose that $S$ is a switch which satisfies $S^2=1$. Then there is a
representation $w_S$ of the flat braid group given by $w_S(\rho_i)=S_i$ and
$w_S(\tau_i)=T_i$. However we can finesse this definition by putting
$w_S(\tau_i)=T'_i$ where $T'$ is any switch which satisfies $T'^2=1$
and $T_1T_2S_1=S_2T_1T_2$. The result will now depend more heavily on
the passage of the representative loop around handles.
\section{Worked examples}

In this section we consider various examples and work out their invariants.
We are very grateful to Andrew Bartholomew who has developed the software to
do the calculations. This can be freely obtained from \cite{B}.

The first example is the projection of the Kishino knot considered in
\cite{FJK}. If we can show that this flat knot is non-trivial
then all possible lifts as virtual knots will {\it a fortiori} be
non-trivial.

\diagram

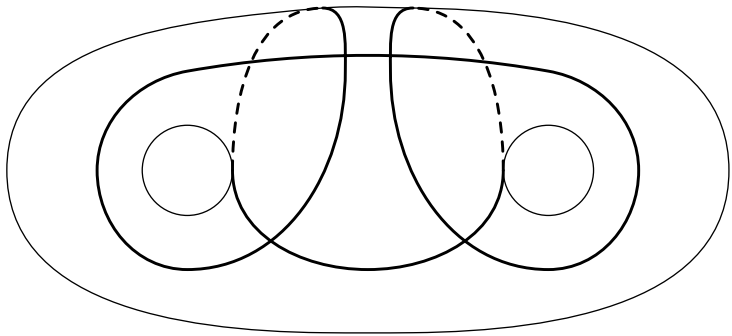
\centerline{\npt The Kishino Flat Knot}

This is the closure of the braid,
$k=\tau_2\sigma_1\sigma_2\sigma_1\tau_2\sigma_1\sigma_2\sigma_1$.
Using the representation of the Weyl algebra given by
$$u=\pmatrix{1&1&0\cr 0&1&1\cr 0&0&1\cr}\quad
v=\pmatrix{y&0&0\cr 1&y&0\cr 0&2&y\cr}$$
with underlying ring $\Z_3[y]$ and using the software
developed by Bartholomew we find
$$\Delta_0=0,\quad \Delta_1=2+2y.$$
Since $\Delta_1\ne1$, the Kishino flat knot is indeed non-trivial.

The second example consists of all possible flat knots which are the
closures of braids in the flat braid group $FB_2$. Clearly we need
only consider closures of
$r_n=\tau_1\sigma_1\tau_1\sigma_1\cdots\tau_1\sigma_1$ where $2n$ is
the number of multiplicands. Let $L_n$ denote this closure.

For example $L_2$, the closure of $r_2=\tau_1\sigma_1\tau_1\sigma_1$ is
illustrated below.
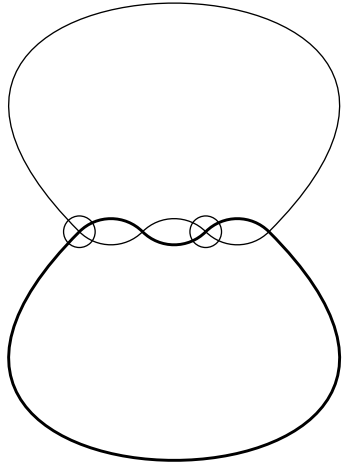
\diagram
We will use the representation of the Weyl algebra given by
$$u=\pmatrix{x&1\cr 0&x\cr}\quad v=\pmatrix{1&0\cr 1&1\cr}$$
with underlying ring $\Z_2[x]$. 

\theorem{\sl With the above representation the invariant of  $L_n$
is $\Delta_0=x^{2n}+x^{-2n}$.}

{\bf Proof} We will look for the eigenvalues of $ST$. This leads to the
equations
$$B\x+A\y=\lambda\x\quad D\x+C\y=\lambda\y$$
where
$$\matrix{
A&=\pmatrix{1/x&1/x^2\cr 1/x&(1+x)/x^2\cr}& 
B&=\pmatrix{x&1\cr 0&x\cr}\cr 
C&=\pmatrix{(1+x+x^2+x^4)/x^5&(1+x^2+x^4)/x^6\cr 1/x^4&(1+x+x^2+x^4)/x^5\cr}&
D&=\pmatrix{(1+x)/x^2&1/x^2\cr 1/x&1/x\cr}\cr}
$$
Eliminating $\y=A^{-1}(\lambda-B)\x$ we see that the matrix
$$D+(C-\lambda)A^{-1}(\lambda-B)$$
is singular. Taking determinants we get the following equation in $\lambda$,
$$\lambda^4+\lambda^2(x^2+x^{-2})+1=0.$$
The roots are $\lambda=x$(twice) and $\lambda=x^{-1}$(twice).
By working over an algebraically closed field extension if necessary we
can write
$$ST=P^{-1}UP$$
where $U$ is upper triangular and has $x,x,x^{-1},x^{-1}$ down the diagonal.
Then $(ST)^n=P^{-1}U^nP$ where $U^n$ is upper triangular and has $x^n,x^n,
x^{-n},x^{-n}$ down the diagonal. So $\Delta_0=\det((ST)^n-1)=x^{2n}+x^{-2n}$
as required. \qed

The third example we call the {\bf whorl}, $W_n$. It is the closure of
the flat braid
$$w_n=\tau_1\tau_2\cdots\tau_n\tau_{n-1}\cdots\tau_2\sigma_1
\sigma_2\cdots\sigma_n.$$
The whorl has genus 2.

\diagram
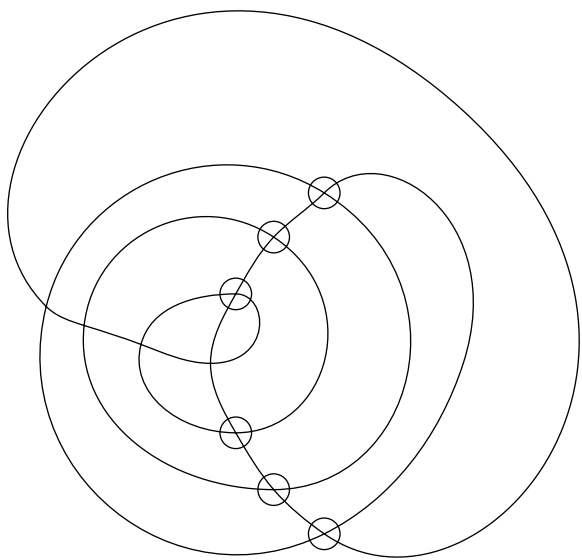
\centerline{\npt The Whorl, $W_4$}

The values of $\Delta_0$ for $n=3, 4$ are
$$\Delta_0 = (x^{10}+x^4+x^2+1)/x^{10}, \quad
(x^{14}+x^{10}+x^8+x^6+x^2+1)/x^{12},$$
respectively. We have been unable to find a general pattern.

\mathsurround=2pt

\centerline{\bf References}

\refe

\cite{As} Helmer Aslaksen, Quaternionic Determinants, Math. Intel. Vol 18 no.
3 (1996)

\cite{B} http://www.layer8.co.uk/maths/braids/

\cite{BF} A. Bartholomew and Roger Fenn. Quaternionic Invariants of Virtual
Knots and Links, to appear in J. Knot Th. Ramifications.  Preprint
available from\nl
http://www.maths.sussex.ac.uk////Staff/RAF/Maths/Current/Andy/

\cite{BuF} S. Budden and Roger Fenn. The equation, $[B,(A-1)(A,B)]=0$
and virtual knots and links, Fund Math 184 (2004) pp 19-29.

\cite{Co} P. M. Cohn. Algebra vol 3 Wiley 1991.

\cite{DP} R. Diaz and E. Pariguan, Symmetric quantum Weyl algebras,
preprint arX QA/0311128 v4 11 Jun 2004

\cite{F} Roger Fenn, Quaternion Algebras and Invariants of Virtual Knots
and Links, preprint

\cite{FJK} R. Fenn, M. Jordan, L. Kauffman,  Biquandles and
Virtual Links, Topology and its Applications, 145 (2004) 157-175

\cite{H} J. Humphreys, Reflection Groups and Coxeter Groups, CUP
(1990)

\cite{Kad} T. Kadokami, Some numerical invariants of flat virtual links are
always realized by reduced diagrams,  to appear in J. Knot Th. Ramifications.

\cite{K} L.Kauffman. Virtual Knot Theory, European J. Comb. Vol 20,
663-690, (1999)

\cite{L} T. Y. Lam. The Algebraic Theory of Quadratic Forms, Benjamin (1973)

\cite{S} D. Stanovsky On axioms of biquandles, to appear in J. Knot Th. 
Ramifications.

\cite{Sa} J. Sawollek, On Alexander-Conway polynomials for virtual knots and links,
preprint. http://citeseer.ifi.unizh.ch/sawollek01alexanderconway.html  

\cite{Tu}  V.Turaev, Virtual strings. Ann. Inst. Fourier 54 (2004),  no. 7,
2455--2525.

\bye